\documentclass[5pt]{article}
\usepackage{amsmath}
\usepackage[latin1]{inputenc}
\usepackage[font=small,labelfont=bf,labelsep=none]{caption}
\usepackage{amsmath,amsthm,amssymb}
\usepackage{graphicx}
\usepackage{color}

\renewcommand{\theequation}{\thesection.\arabic{equation}}

\newtheorem{lemma}{Lemma}[section]
\newtheorem{thm}{Theorem} [section]

\newtheorem{exmp}{Example} [section]
\newtheorem{cor}{Corollary}[section]
\newtheorem{rem}{Remark}[section]
\title{ Two Generalizations for Quadratic Residue Codes over Finite Fields}
\author{{
\quad Qunying Liao\footnote{ \small Q.Y.Liao, College of Mathematical Science, Sichuan Normal University,
  \small Chengdu, Sichuan, China(Email: qunyingliao@sicnu.edu.cn)} \quad Yuanbo Liu \footnote{Y.B.Liu, College of Mathematical Science, Sichuan Normal University,
 \small Chengdu, Sichuan, China(Email: 361717444@qq.com)} }\ }
\date{}
\begin{document}
\baselineskip15pt \maketitle
\renewcommand{\theequation}{\arabic{section}.\arabic{equation}}
\catcode`@=11 \@addtoreset{equation}{section} \catcode`@=12
\begin{abstract}

It's well known that the quadratic residue code over finite fields is an interesting class of cyclic codes for its higher minimum distance. Let $g$ be a positive integer and $p,p_{1},\ldots, p_{g}$ be distinct odd primes, the present paper generalizes the constructions for the quadratic residue code with length $p$ to be the length $n=p_{1}\cdots p_{g}$, and to be the case $m$-th residue codes with length $p$ over finite fields, where $m\geq 2$ is a positive integer. Furthermore, a criterion for that these codes are self-orthogonal or complementary dual is obtained, and then the corresponding counting formula are given. In particular, the minimum distance of all 24 quaternary quadratic residue codes $[15,8]$ are determined.


\end{abstract}

\vskip6mm \section{Introduction and Backgrounds}

Error-correcting codes, especially linear codes, are widely applied for many fields, such as for networking, satellite communication, compact disks, and so on.
Let $q$ be a power of the prime and $\mathbb{F}_{q}$ denote the finite field of $q$ elements. Refer to \cite{M.S}, a $q$-ary linear code $\mathcal{C}$ of length $n$ is a $k$-dimension subspace of $\mathbb{F}_{q}^{n}$. The vector $\mathbf{c}=(c_{0},\ldots,c_{n-1})\in \mathcal{C}$ is called a codeword and the weight of a codeword $\mathbf{c}$ is the number of nonzero $c_{i}\neq 0 \ (0\leq i\leq n-1)$. The minimum (Hamming) distance $d$ of a linear code $\mathcal{C}$ is the minimum weight value of all nonzero codewords in $\mathcal{C}$. For applications, the value $\frac{k}{n}$ represents the translation efficiency of $\mathcal{C}$, and the minimum distance $d$ decides the capability of error-checking and correcting. Therefore a $q$-ary linear code $\mathcal{C}$ of length $n$, dimension $k$ and minimum distance $d$ is denoted to be $\mathcal{C}=[n,k,d]_q$ or $\mathcal{C}=[n,k,d]$ when the finite field is fixed.

It's well-known that dual codes are also important for applications. The dual code of $\mathcal{C}$ is defined to be
$$
\mathcal{C}^\bot =\{b\in \mathbb{F}_{q}^{n}| <b,c>=0, \forall c\in \mathcal{C}\},
$$
where $<b,c>$ denotes the inner product of two vectors $b$ and $c$. It is easy to see that for a $q$-ary linear code $\mathcal{C}=[n,k,d]$, the dual code $\mathcal{C}^{\bot}$ is also a $q$-ary linear code with length $n$ and dimension $n-k$. Furthermore, a $q$-ary linear code $\mathcal{C}$ is called a self-orthogonal code if $\mathcal{C}\cap\mathcal{C}^{\bot}=\mathcal{C}^{\bot}$, and is called a linear code with complementary dual code, in short LCD, if $\mathcal{C}\cap\mathcal{C}^{\bot}=\{0\}$, equivalently, the direct sum subspace $\mathcal{C}\oplus\mathcal{C}^{\bot}=\mathbb{F}_{q}$. In recent years, there are many good results for self-orthogonal codes or LCD over finite fields $\cite{J.L,L.D}$.

One of the main problems in coding theory is to find those codes with the efficiency $\frac{k}{n}$ and the capability $d$ both large \cite{M.S}. While these parameters are restricted to each other, namely, which is showed by some bounds such as GV Bound, Hamming Bound, Singleton Bound, and so on. The well-known Singleton Bound says that $d\leq n+1-k$ for any code $\mathcal{C}=[n,k,d]_q$ and $\mathcal{C}$ is called a maximum distance separable (MDS) code when the equality is true. Especially, as a class of linear codes, cyclic codes over finite fields have wide applications in storage and communication systems due to their efficient encoding and decoding algorithms \cite{C,F,P,R}. It has been studied for decades and a lot of progress has been made \cite{H. P}. A $q$-ary linear code $\mathcal{C}=[n,k]$ is cyclic if any $(c_{0},c_{1},\ldots,c_{n-1})\in \mathcal{C}$ implies $(c_{n-1},c_{0},\ldots,c_{n-2})\in \mathcal{C}.$ By the one-to-one correspondence
$$\varphi : \mathcal{C} \rightarrow R:=\mathbb{F}_{q}[x]/(x^{n}-1)$$
$$(c_{0},c_{1},\ldots,c_{n-1})\longmapsto c_{0}+c_{1}x+ \cdots +c_{n-1}x^{n-1}.$$
Evidently, any $q$-ary cyclic code $\mathcal{C}$ with length $n$ is equivalent to be an ideal of the ring $R$. Note that $R$ is a principal ideal ring, then there exists a unique monic polynomial $g(x)$ with the least degree such that $\varphi(\mathcal{C})=<g(x)>$, where $g(x)$ is a monic factor of $x^n-1$ with degree $n-k$. In general, $g(x)$ is called the generator polynomial of $\mathcal{C}$ and $h(x):=\frac{x^{n}-1}{g(x)}$ is called the parity-check polynomial of $\mathcal{C}.$ For a polynomial $f(x)=\sum_{i=0}^{k}f_{i}x^{i}\in \mathbb{F}_{q}[x]$ with degree $k$ and $f_{0}\neq 0,$ the reciprocal polynomial of $f(x)$ is defined to be $f^{*}(x):=f_{0}^{-1}x^{k}\sum _{i=0}^{k}f_{i}(\frac{1}{x})^{i}.$ Then $f(\alpha)=0$ if and only if $f^{*}(\alpha^{-1})=0,$ and $f(x)$ is called self-reciprocal if $f(x)=f_{k}f^{*}(x).$ It is easy to see that for a cyclic code $\mathcal{C}=<g(x)>$, the dual $\mathcal{C}^{\bot}=<h^*(x)>$.

\textbf{Proposition 1.1}$^{\cite{L.D}}$\quad Let $n$ be a positive integer and $q$ be the power of a prime with $\gcd(n,q)=1$. Suppose that $\mathcal{C}=<g(x)>$ is a $q$-ary cyclic code with length $n$, and $\Omega_q$ is the algebraic closure of $\mathbb{F}_q$. Then the following statements are equivalent.

(1)\quad $\mathcal{C}$ is an LCD code.

(2)\quad $g(x)$ is a self-reciprocal polynomial in $\mathbb{F}_q[x]$.

(3)\quad For any $\alpha\in\Omega_q$ with $\alpha^n=1, g(\alpha)=0$ means that $g(\alpha^{-1})=0$.

\textbf{Proposition 1.2}\quad Let $n$ be a positive integer, $q$ be the power of a prime with $\gcd(n,q)=1$, and $\Omega_q$ be the algebraic closure of $\mathbb{F}_q$. Suppose that $\mathcal{C}=<g(x)>$ is a $q$-ary cyclic code with length $n$, and $h^{*}(x)$ is the reciprocal polynomial of the parity-checking polynomial $h(x)$ of $\mathcal{C}$. Then the following statements are equivalent.

(1)\quad $\mathcal{C}$ is a self-orthogonal code.

(2)\quad $g(x)$ is a factor of $h^{*}(x)$.

(3)\quad  For any $\alpha\in\Omega_q$ with $\alpha^n=1, g(\alpha)=0$ means that $g(\alpha^{-1})\neq 0$.

\textbf{Proof}\quad  The equivalence of (1) and (2) can be seen in \cite{P.Q}.

Next, we prove that (2) is equivalent to (3).

If (2) is true, then for any $\alpha\in\Omega_q, g(\alpha)=0$ means that $h^{*}(\alpha)=0$, and then $h(\alpha^{-1})=0$. While from assumptions we know that $x^n-1=g(x)h(x)$ has no repeated roots in $\Omega_q$, therefore $g(\alpha^{-1})\neq 0$, which means that (3) is true.

Now suppose that (3) is true, then by $\gcd(n,q)=1$ and $x^n-1=g(x)h(x)$, it is enough to show that for any $\alpha\in\Omega_q, g(\alpha)=0$ means that $h^{*}(\alpha)=0$. In fact, from $g(\alpha)=0$ and $x^n-1=g(x)h(x)$ we know that $\alpha^n=1$ and then $\alpha^{-n}=1$, hence $g(\alpha^{-1})h(\alpha^{-1})=0$. But from (3) and $g(\alpha)=0$ we have $g(\alpha^{-1})\neq 0$, thus we must have $h(\alpha^{-1})=0$, which means that $h^{*}(\alpha)=0$, namely, (2) is true.

This completes the proof of Proposition 1.2.  \hfill$\Box$\\

Based on the generator of a cyclic code, one immediately gets two parameters $n$ and $k$, another important parameter $d$ can be obtained from the parity-checking polynomial. In fact, if the parity-checking polynomial $h(x)=a_{0}+a_{1}x+\cdots +a_{k}x^{k}$ of a cyclic code $\mathcal{C}=[n,k]_q$, one can get a parity-checking matrix of $\mathcal{C}$ as follows:
\begin{equation}
 H=\left(
 \begin{array}{ccccc}
 a_{0} & a_{1} & a_{2} & \cdots & 0  \\
  0 & a_{0} & a_{1} & \cdots & 0  \\
  0 & 0 & a_{0} & \cdots & 0  \\
  \vdots & \vdots & \vdots & \vdots & \vdots \\

  0 & 0 & 0 & \cdots  & a_{k}  \\
  \end{array}
 \right)_{(n-k) \times n}.
\end{equation}

\begin{lemma}$^{\cite{H. P}}$
{ If $H$ is the parity check matrix of a $q$-ary linear code $\mathcal{C}$, then $\mathcal{C}$ has minimum distance $d$ if and only if any $d-1$ columns of $H$ are linearly $\mathbb{F}_{q}$-independent and some $d$ columns are linearly  $\mathbb{F}_{q}$-dependent.
}
\end{lemma}

On the other hand, the quadratic residue code over finite fields, introduced by Prange in 1958 \cite{P}, is a type of cyclic codes with a prime length. Since it's code rate is greater than or equal to $\frac{1}{2}$ and has higher minimum distance \cite{C.R.T}, several decoding techniques for classic quadratic residue codes have been developed in the past decades. Correspondingly, some generalized constructions are obtained. For examples, in 2006, based on properties for cyclotomic fields, a class of quadratic residue codes with a prime length are constructed \cite{E.I}. In 2009, Charters\cite{C.P} generalized binary quadratic residue codes with a prime length over finite fields to be the case $p$-th residue codes, where $p$ is a prime, and then obtained a lower bound for the codeword weight of these codes and their dual codes. In 2012, for two distinct odd primes $p_1$ and $p_2$, Ding \cite{D} introduced a cyclotomic construction for cyclic codes with length $p_{1}p_{2}$ over a binary field, ternary field or quaternary field, respectively, and gave a lower bound, which is the best bound till now, of the minimum odd-like weight for these codes. In 2017, Li., et al. \cite{L.D}, constructed several classes of LCD cyclic codes over finite fields.

The present paper continues the study for quadratic residue codes over finite fields, generalizes Ding's constructions for quadratic residue codes with a prime length in 2012 to be the length $n=p_{1}\cdots p_{g}$, where $p_{1},\ldots,p_{g}$ are distinct odd primes, and then obtains a criterion and the computing formula for these (self-orthogonal or LCD) quadratic residue codes. In particular, we completely determine the minimum distance of all 24 quadratic residue code $[15,8]_4$ listed in Table 1.
\begin{table}[!h]
\centering
\caption{}
\begin{tabular}{|c|c|}
 \hline
   Number of $[15,8]_{4}$ codes & Minimum Distance \\
  $16$                    &6 \\
  $4$                     & 4\\
  $4$                     &3 \\

  \hline
\end{tabular}
\end{table}
Meanwhile, in section 3, we generalize the constructions of Charters and Ding for quadratic residue codes with length $p$ to be the $m$-th residue code over finite fields, where $p$ is an odd prime and $m\geq 2$ is a positive integer with $m|p-1$, and then obtain a criterion and the corresponding computing formula for these (LCD or self-orthogonal) codes.

\section{QR Codes with Odd Length over $\mathbb{F}_{q}$}

\subsection{The Constructions for Quadratic Residue Codes over $\mathbb{F}_{q}$}

For convenience, throughout this section, suppose that $g$ is a positive integer, $p_{1},\ldots,p_g$ are distinct odd primes, $n=p_{1}\cdots p_{g},$  and $q$ is a power of the prime $p$ with $\gcd(n,p)=1$. We also define
$$
N=ord_{n}(q)=lcm(ord_{p_{1}}(q),\ldots,ord_{p_{g}}(q)),
$$
 and $\theta=\alpha^{\frac{q^{N}-1}{n}}$, where $\alpha$ is a primitive element of $\mathbb{F}_{q^{N}}^{*}$. Hence
$\theta$ is a fixed $n$-th primitive unit root  in $\mathbb{F}_{q^{N}}^{*}$.

First, we give a partition of the set $M=\{1,2,\ldots, n\}$. Set
$$
A_0=\{1\}, M_{1}=\{i \,|\,\gcd (i,n)=1,1\leq i \leq n-1\},
$$
$$
A_r=\{p_{j_{1}}p_{j_{2}} \cdots p_{j_{r}} ,1\leq j_{1}<\cdots <j_{r}\leq g-1\}(r=1,2, \ldots , g-1),$$
and for any $i_r\in A_r(1\leq r\leq g-1)$, denote
$$
M_{i_r}=\{j\,|\,\gcd(j,n)=i_r, 1\leq j\leq n-1\}.
$$
Thus we can get the partition
$$
M=\{n\}\cup M_{1}\cup_{r=1}^{g-1}\cup_{i_r\in A_r}M_{i_r}=\{n\}\cup_{r=0}^{g-1}\cup_{i_r\in A_r}M_{i_r},
$$
and then the factorization
$$
x^{n}-1=(x-1)\prod_{r=0}^{g-1}\prod_{i_r\in A_r}\prod_{t\in M_{i_r}}(x-\theta^t).
$$

Secondly, for any $r=0,1,\ldots,g-1$, set
$$
M^*_{i_r}=\left\{Q_{i_r}| Q_{i_r}|\frac{n}{i_r}, 1<Q_{i_r}\leq n\right\}.$$
And for given $Q_{i_r}\in M^*_{i_r}(i_r\in A_r)$, denote
$$
M_{i_r,Q_{i_r}}^{1}=\left\{j\,|\, \left(\frac{j}{Q_{i_r}}\right)=1, j\in M_{i_r}\right\},\quad  M_{i_r,Q_{i_r}}^{-1}=\left\{j\,|\, \left(\frac{j}{Q_{i_r}}\right)=-1, j\in M_{i_r}\right\},
$$
then we have the disjoint unions
$$
M_{i_r}=M_{i_r,Q_{i_r}}^{1}\cup M_{i_r,Q_{i_r}}^{-1}.
$$

Finally, for any  $\epsilon_{Q_{i_r}}\in \{1,-1\}$, denote
$$
 F^{\epsilon_{Q_{i_r}}}_{i_r,Q_{i_r}}(x)=\prod _{t\in M^{\epsilon_{Q_{i_r}}}_{i_r,Q_{i_r}}}(x-\theta ^{t}).
$$

Now, we can get the following generalized definition for quadratic residue codes with length $n$ over finite fields.

\textbf{Definition 2.1}\quad  For any $r=0,1,\ldots,g-1$, $\epsilon_{Q_{i_r}}\in \{1,-1\}$ and fixed $Q_{i_r}\in M^*_{i_r}$, the cyclic code
$$
\mathcal{C}=<g(x)>=<\prod_{r=0}^{g-1}\prod_{i_r\in A_r}F^{\epsilon_{Q_{i_r}}}_{i_r,Q_{i_r}}(x)>
$$
is called a quadratic residue (QR) code with length $n$.

It is easy to see that the QR code in Definition 2.1 has the largest dimension $\frac{n+1}{2}$, which is just the classical quadratic residue code for $g=1$, and coincides with Ding's cyclotomic construction for $g=2$ \cite{D}. In particularly, for $Q_{1}\in\{p_{1},p_{2},n\}, Q_{p_1}=p_{2}$ and $Q_{p_2}=p_{1}$, one can get eight cyclic codes, which coincides with Ding's construction. And for $Q_{1}=n$, it is just the first construction for eight cyclic codes in Propositions 2, 5-7  given by Ding  \cite{D}. When $Q_{1}=p_1$ or $p_2$, it is exactly the second construction for eight cyclic codes in Propositions 8, 11-13 or the third construction for eight cyclic codes in Propositions 14, 17-19 \cite{D}, respectively. While for the case $g\geq 3$, our construction is new and different from Ding's construction in \cite{D}. In section 2.2, we give a sufficient and necessary condition for the existence of LCD or self-orthogonal QR codes over finite fields (Theorem 2.2), and then obtain the counting formula for QR codes with length $n$, which is the generalization of the corresponding results obtained by Ding in 2012. Furthermore, we obtain the explicit number of LCD or self-orthogonal QR codes with length $n$ over finite fields (Theorem 2.3), respectively.

\subsection{Main Results}

In this section, we prove that the quadratic residue code $\mathcal{C}$ in Definition 2.1 is a $q$-ary code when the Legendre symbol $\left(\frac{q}{p_{i}}\right)=1$ for any $i=1,\ldots,g$. It is enough to show that the generator polynomial $g(x)\in\mathbb{F}_{q}[x]$, equivalently,
$$
F^{\epsilon_{Q_{i_r}}}_{i_r,Q_{i_r}}(x)\in \mathbb{F}_{q}[x],
$$
where $ \epsilon_{Q_{i_{r}}}\in \{1,-1\}$ for any $r=0,1,\ldots,g-1, i_r\in A_r$ and  $Q_{i_r}\in M^*_{i_r}.$  The fact can be deduced from the following theorem directly.

\begin{thm}

{The assumptions are as above. For each $ \epsilon_{Q_{i_{r}}}\in \{1,-1\}$,  we have
 $$
 F^{\epsilon_{Q_{i_r}}}_{i_r,Q_{i_r}}(x)\in \mathbb{F}_{q}[x]/(x^{n}-1)\Longleftrightarrow
 q M_{i_r,Q_{i_r}}^{\epsilon_{Q_{i_{r}}}}= M_{i_r,Q_{i_r}}^{\epsilon_{Q_{i_{r}}}} \Longleftrightarrow \left(\frac{q}{p_{i}}\right)=1(1\leq i\leq g),
 $$
 where $j\in q M_{i_r,Q_{i_r}}^{\epsilon_{Q_{i_{r}}}}$ is taken to be the minimum non-negative residue modulo $n$.}
\end{thm}

From Theorem 2.1 and properties for the Legendre or Jacobi symbols, we immediately have the following results.
\begin{cor}
{For a quadratic residue code $\mathcal{C}$ in Definition 2.1,

(1)\quad $\mathcal{C}$ is a quaternary quadratic residue code;

(2)\quad  $\mathcal{C}$ is a binary or a ternary quadratic residue code, if and only if for any $i=1,2,\ldots,g,$
 $$
 p_{i} \equiv \pm 1 (\mbox{mod}\ 8) \ \mbox{or} \quad  p_{i} \equiv \pm 1(\mbox{mod}\ 12),
 $$
 respectively.}
\end{cor}

\begin{thm}

{(1)\quad The number of all QR codes with length $n$ is $N=\prod_{r=0}^{g-1}(2\times(2^{g-r}-1))^{\dbinom{g}{r}}$.

(2)\quad There exists a LCD QR code if and only if $p_{i}\equiv 1(\mbox{mod}\ 4)$ for any $i=1,\ldots,g$.

(3)\quad  There exists a self-orthogonal QR code if and only if $p_{i}\equiv 3(\mbox{mod}\ 4)$ for any $i=1,\ldots,g$.}
\end{thm}

\begin{rem}

{By taking $g=1$, from (1) of Theorem 2.2 there are exactly $2$ QR codes with a prime length, which is just the corresponding result given by Prange in 1958.

(2)\quad By taking $g=2$, from (1) of Theorem 2.2 there are exactly $24$ QR codes with length $n=g_1g_2$, which is just the corresponding result given by Ding in 2012.}
\end{rem}

\begin{thm}

{(1)\quad The number of LCD QR codes with length $n$ equals to $0$ or $N$.

(2)\quad The number of self-orthogonal QR codes with length $n$ equals to $0$ or
$$
T=2^{g2^{g-1}}.
$$}
\end{thm}

\subsection{The Proofs For Main Results}

\textbf{The Proof for Theorem 2.1}.

For any $\epsilon_{Q_{i_{r}}}\in \{1,-1\}(r=0,1,\ldots,g-1, i_r\in A_r$ and  $Q_{i_r}\in M^*_{i_r})$,  from the definitions we know that
$$F^{\epsilon_{Q_{i_r}}}_{i_r,Q_{i_r}}(x)=\prod _{t\in M^{\epsilon_{Q_{i_r}}}_{i_r,Q_{i_r}}}(x-\theta ^{t}), \ (F^{\epsilon_{Q_{i_r}}}_{i_r,Q_{i_r}}(x))^{q}=\prod _{t\in M^{\epsilon_{Q_{i_r}}}_{i_r,Q_{i_r}}}(x^{q}-\theta ^{tq})$$
and
$$
F_{rj}^{\epsilon_{rj}}\in  \mathbb{F}_{q}[x]/(x^{n}-1)\Longleftrightarrow (F^{\epsilon_{Q_{i_r}}}_{i_r,Q_{i_r}}(x))^{q}=F^{\epsilon_{Q_{i_r}}}_{i_r,Q_{i_r}}(x^{q}),
$$
i.e.,
$$
F_{rj}^{\epsilon_{rj}}\in  \mathbb{F}_{q}[x]/(x^{n}-1) \Longleftrightarrow \prod _{t\in M^{\epsilon_{Q_{i_r}}}_{i_r,Q_{i_r}}}(x^{q}-\theta ^{tq})=\prod _{t\in M^{\epsilon_{Q_{i_r}}}_{i_r,Q_{i_r}}}(x^{q}-\theta ^{t}).
$$
Note that $\theta$ is an $n$-th primitive root of unity in $\mathbb{F}_{q^{N}}^{*}$, hence
$$
\prod _{t\in M^{\epsilon_{Q_{i_r}}}_{i_r,Q_{i_r}}}(x^{q}-\theta ^{tq})=\prod _{t\in M^{\epsilon_{Q_{i_r}}}_{i_r,Q_{i_r}}}(x^{q}-\theta ^{t}) \Longleftrightarrow qM^{\epsilon_{Q_{i_r}}}_{i_r,Q_{i_r}}= M^{\epsilon_{Q_{i_r}}}_{i_r,Q_{i_r}},
$$
which means that
$$
 F^{\epsilon_{Q_{i_r}}}_{i_r,Q_{i_r}}(x)\in \mathbb{F}_{q}[x]/(x^{n}-1)\Longleftrightarrow
 q M_{i_r,Q_{i_r}}^{\epsilon_{Q_{i_{r}}}}= M_{i_r,Q_{i_r}}^{\epsilon_{Q_{i_{r}}}}. \eqno(3.1)
$$

Furthermore, if $q M_{i_r,Q_{i_r}}^{\epsilon_{Q_{i_{r}}}}=M_{i_r,Q_{i_r}}^{\epsilon_{Q_{i_{r}}}}$, especially, $qM_{i_{g-1},Q_{i_{g-1}}}^{1}=M_{i_{g-1},Q_{i_{g-1}}}^{1}$ for each
$$
i_{g-1}\in A_{g-1}=\{\frac{n}{p_{1}},\frac{n}{p_{2}}, \ldots ,\frac{n}{p_{g}}\}
$$
and $Q_{i_{g-1}}\in M^*_{i_{g-1}}$. From the assumptions we know that $Q_{i_{g-1}}=p_{i}$ for each $i_{g-1}=\frac{n}{p_{i}}(i=1,2,\ldots, g)$, and
$$
M_{i_{g-1},Q_{i_{g-1}}}^{1}=\left\{ j\, |\, \gcd(j,n)=i_{g-1}=\frac{n}{p_{i}},1\leq j \leq n-1,\left(\frac{j}{p_i}\right)=1\right\}.
$$
And so $qM_{i_{g-1},Q_{i_{g-1}}}^{1}=M_{i_{g-1},Q_{i_{g-1}}}^{1}$, which means that for any $j\in M_{i_{g-1},Q_{i_{g-1}}}^{1}, \left(\frac{j}{p_{i}}\right)=1\Longleftrightarrow\left(\frac{jq}{p_{i}}\right)=1,$ i.e.,
$$
\left(\frac{q}{p_{i}}\right)=1\  (1\leq i\leq g).\eqno(3.2)
$$

On the other hand, if (3.2) is true, then from $Q_{i_r}|\frac{n}{i_r}\ (0\leq r\leq g-1,i_r\in A_r)$,
we have
$$
 \left(\frac{q}{Q_{i_r}}\right)=1\ (0\leq r\leq g-1,i_r\in A_r).
$$
Now from the assumptions that $\gcd(n,q)=1$ and for $\epsilon_{Q_{i_r}}\in \{1,-1\}$,
$$
M_{i_r,Q_{i_r}}^{\epsilon_{Q_{i_r}}}=\left\{ j \,| \,\gcd(j,n)=i_r, i_r\in A_{r} ,1\leq j \leq n-1,\left(\frac{i}{Q_{i_r}}\right)=\epsilon_{Q_{i_r}}\right\},
$$
we can get $\gcd(jq,n)=i_r$ and $\left(\frac{jq}{Q_{i_r}}\right)=\epsilon_{Q_{i_r}}$ for any $j\in M_{i_r,Q_{i_r}}^{\epsilon_{Q_{i_r}}},$ namely,
$$
q M_{i_r,Q_{i_r}}^{\epsilon_{Q_{i_r}}}=M_{i_r,Q_{i_r}}^{\epsilon_{Q_{i_r}}} (r=0,1,\ldots,g-1, i_r\in A_r , Q_{i_r}\in M^*_{i_r}).\eqno(3.3)
$$

Hence, from (3.1)-(3.3) we complete the proof of Theorem 2.1. \hfill$\Box$\\

\textbf{Proof for Theorem 2.2}.

(1)\quad By the definitions of $A_r$ and $M^*_{i_r}$, we have
$$
|A_r|=\dbinom{g}{r},\quad |M^*_{i_r}|=2^{g-r}-1.
$$

For $r=0,1,\ldots,g-1$ and $Q_{i_r}\in M^*_{i_r}, F^{\epsilon_{Q_{i_r}}}_{i_r,Q_{i_r}}$ has exactly two choices. And then there are $2|M^{*}|=2(2^{2^{g-r}}-1)$ choices for $F^{\epsilon_{Q_{i_r}}}_{i_r,Q_{i_r}}$.
Hence by Definition 2.1,
 the number of all QR codes with length $n$ is equal to
$$
N=\prod_{r=0}^{g-1}(2\times|M^*_{i_r}| )^{|A_r|}=\prod_{r=0}^{g-1}(2\times(2^{g-r}-1) )^{\dbinom{g}{r}}.
$$
This completes the proof of (1).

(2)\quad By Proposition 1.2, $\mathcal{C}$ is LCD if and only if $\beta^{-1}$ is also a root of $g(x)$ when $\beta$ is a root of $g(x)$, which means that $\beta$ and $\beta^{-1}$ are both roots of $F^{\epsilon_{Q_{i_r}}}_{i_r,Q_{i_r}}$ for some $\epsilon_{Q_{i_r}}\in \{1,-1\}$. Note that $\beta$ is an $n$-th root of unity and $\theta\in \mathbb{F}_{q}$ is an $n$-th primitive root of unity. So $\beta=\theta^i$ for some $i=1,\ldots,n-1$, and then $i$ and $-i$ are both in the same set $M_{i_r,Q_{i_r}}^{\epsilon_{Q_{i_r}}}$ for some $Q_{i_r}\in M^*_{i_r}$. Therefore there exists a LCD QR code with length $n$ if and only if for any $r=0,1,\ldots,g-1, i_r\in A_r$ and some $Q_{i_r}\in M^*_{i_r}$, we have $\left(\frac{-1}{Q_{i_r}}\right)=1$.

Now suppose that there exists a LCD QR code $\mathcal{C}$, then for any $i_{g-1}\in A_r,$ we have $\left(\frac{-1}{Q_{i_{g-1}}}\right)=1$ for some $Q_{i_{g-1}}\in M^*_{i_{g-1}}$. From the definition we know that
$$
A_{g-1}=\left\{\frac{n}{p_i}\,|\, 1\leq i\leq g\right\}.
$$
This means there are exactly $g$ sets $M^*_{i_{g-1}}=\{p_j\}(j=1,\ldots,g)$. Thus we can get
$$
\left(\frac{-1}{p_1}\right)= \cdots =\left(\frac{-1}{p_g}\right)=1,
$$
which is true if and only if $p_{i}\equiv 1\ \ (\mbox{mod} \ 4)$ for any $i=1,\ldots,g$.

On the other hand, by the definition of $M^*_{i_r}$ we know that $Q_{i_r}\in M^*_{i_r}$ is a product of some $p_i(1\leq i\leq g)$, thus by properties of the Jacobi symbol and the assumption that
$p_{i}\equiv 1\ (\mbox{mod} \ 4)$ for any $i=1,\ldots,g$, we can get $\left(\frac{-1}{Q_{i_r}}\right)=1$ for any $i_r\in A_r$ and $Q_{i_r}\in M^*_{i_r}.$ And so there exists a LCD QR code with length $n$.

 This completes the proof of (2).

(3)\quad In the same proof as that for (2), there exists a self-orthogonal QR code if and only if for some $i=1,\ldots,n-1,$ one of $i$ and $-i$ is included in $M_{i_r,Q_{i_r}}^1$ and the other is included in $M_{i_r,Q_{i_r}}^{-1}$ for some $r\in \{0,\ldots,g-1\}$. Therefore there exists a self-orthogonal QR code if and only if for any $r=1,\ldots,g-1, i_r\in A_r$ and some $Q_{i_r}\in M^*_{i_r}$, we have $\left(\frac{-1}{Q_{i_r}}\right)=-1$. Now by the definition of $Q_{i_r}$ and properties of the Jacobi symbol, the statement is true if and only if for any $i=1,\ldots,g, p_{i}\equiv 3\  (\mbox{mod} \ 4)$, and for any $r=0,\ldots,g-1, Q_{i_r}=\prod^k_{s=1}p_{j_s}$ for some odd integer $k$. This completes the proof of (3). \hfill$\Box$\\

\textbf{Proof for Theorem 2.3}.

(1)\quad By the proof of Theorem 2.2 (2), there exists a LCD QR code if and only if any QR code is LCD. Hence by (1) of Theorem 2.2, (1) is immediate.

(2)\quad  Suppose that there exists a self-orthogonal QR code with length $n$, then from the proof of Theorem 2.2 (3), we know that the QR code
$$
\mathcal{C}=<g(x)>=<\prod_{r=0}^{g-1}\prod_{i_r\in A_r}F^{\epsilon_{Q_{i_r}}}_{i_r,Q_{i_r}}(x)>
$$
is self-orthogonal if and only if $\left(\frac{-1}{Q_{i_r}}\right)=-1$ for any $r=0,\ldots,g-1, i_{r}\in A_r$, which means that $\left(\frac{-1}{p_{i}}\right)=-1(i=1,2,\cdots ,g)$ and $ Q_{i_r}=\prod^k_{s=1}p_{j_s}$ with some odd $k$.

Now, for any $r=0,1,\ldots,g-1$ and $i_r\in A_r$, set
$$
H_{i_r}=\{Q_{i_r}\,|\, Q_{i_r}\in M^*_{i_r}\ \mbox{and}\ Q_{i_r}=\prod^k_{s=1}p_{j_s}\ \mbox{for some odd integer}\ k, 1\leq k\leq g\}.
$$
Then for given $i_r\in A_r$, it is easy to show that
$$
|H_{i_r}|=2^{g-r-1}.
$$
Thus from the proof of Theorem 2.2 (1), when $p_i\equiv 3 \ (\mbox{mod}\ 4)$ for any $i=1,\ldots,g$, the number of self-orthogonal QR codes with length $n$ is equal to
$$
T_{0}=\prod^{g-1}_{r=0}(2\cdot|H_{r_i}|)^{|A_r|}=\prod^{g-1}_{r=0}(2\cdot 2^{g-r-1})^{\dbinom{g}{r}}=2^{g2^{g-1}}.
$$
This completes the proof of (2).\hfill$\Box$\\

\subsection{ Some Examples }
In this section, several examples are given to support Theorems 2.1-2.3. We obtain all the quaternary QR codes with length 15, the  binary QR codes with length 161 and the ternary QR codes with length 253 in Example 2.1-2.3, respectively. Based on the extensive calculation, the minimal distance for all 24 QR codes in Example 2.1 are obtained. Furthermore, we also consider QR codes with length three odd primes product in Examples 2.4-2.5, which can test Theorems 2.2-2.3.

\begin{exmp} Consider quaternary QR codes with length $n=15=3\cdot 5$ and dimension $8$. Suppose that $\alpha^{2}+\alpha+1=0$ and $\beta^{2}+\beta+\alpha=0$, then $\alpha$ and $\beta$ are primitive elements of $\mathbb{F}_{2^{2}}$ and $\mathbb{F}_{4^{2}}$, respectively. Then by definitions we have
$$A_0=\{1\},\ A_1=\{3,5\},$$
$$
M_1=\{1,2,4,7,8,11,13,14 \},\quad M_3=\{3,6,9,12\}, \quad M_5=\{5,10\},
$$
$$
M_1^{*}=\{3,5,15\},\quad M_3^{*}=\{5\}, \quad M_5^{*}=\{3\},
$$
and
$$
Q_{1}\in \{3,5,15\},Q_{3}=5 ,Q_{5}=3.
$$

Case 1. By taking $Q_{1}=15,$ we have
$$
M_{1,15}^{1}=\{1,2,4,8\},M_{1,15}^{-1}=\{7,11,13,14 \},M_{3,5}^{1}=\{6,9\},M_{3,5}^{-1}=\{3,12\},M_{5,3}^{1}=\{10\},M_{5,3}^{-1}=\{5\},
$$
and then
$$
x^{15}-1=(x-1)F_{1,15}^{1}(x)F_{1,15}^{-1}(x)F_{3,5}^{1}(x)F_{3,5}^{-1}(x)F_{5,3}^{1}(x)F_{5,3}^{-1}(x),
$$
where
$$
F_{1,15}^{1}(x)=x^{4}+ x+1,F_{1,15}^{-1}(x)=x^{4}+ x^{3}+1,F_{3,5}^{1}(x)=x^{2}+\alpha x+1,F_{3,5}^{-1}(x)=x^{2}+(\alpha +1)x+1,
$$
and
$$
F_{5,3}^{1}(x)=x+\alpha +1,\quad F_{5,3}^{-1}(x)=x+\alpha.
$$

Hence by Definition 2.1 we can get 8 quaternary QR codes with the generator polynomial
$$g(x)=F_{1,15}^{\epsilon _{Q_{1}}}(x)F_{3,5}^{\epsilon _{Q_{3}}}(x)F_{5,3} ^{\epsilon _{Q_{5}}}(x)$$
for $\epsilon _{Q_{i}}\in \{1,-1\}(i=1,3,5).$

Case 2. By taking $Q_{1}=3$, we can get $M_{1,3}^{1}=\{1,4,7,13\}, M_{1,3}^{-1}=\{2,8,11,14\},$ and then
$$
F_{1,3}^{1}(x)=(x^{2}+ x+\alpha)(x^{2}+\alpha x+\alpha), F_{1,3}^{-1}(x)=(x^{2}+ x+\alpha+1)(x^{2}+(\alpha+1)x+\alpha+1).$$

Similarly, by Definition 2.1 we can also get 8 quaternary QR codes with the generator polynomial
$$g(x)=F_{1,3}^{\epsilon _{Q_{1}}}(x)F_{3,5}^{\epsilon _{Q_{3}}}(x)F_{5,3}^{\epsilon _{Q_{5}}}(x)$$
for $\epsilon _{Q_{i}}\in \{1,-1\}(i=1,3,5)$.

Case 3. By taking $Q_{1}=5$, we have $M_{1,5}^{1}=\{1,4,11,14\}, M_{1,5}^{-1}=\{2,7,8,13\},$ and then
$$
F_{1,5}^{1}(x)=(x^{2}+ x+\alpha)(x^{2}+(\alpha+1)x+\alpha+1),F_{1,5}^{-1}(x)=(x^{2}+ x+\alpha+1)(x^{2}+\alpha x+\alpha).
$$
In the same way, by Definition 2.1 we can get 8 quaternary QR codes with the generator polynomial
$$g(x)=F_{1,5}^{\epsilon _{Q_{1}}}(x)F_{3,5}^{\epsilon _{Q_{3}}}(x)F_{5,3}^{\epsilon _{Q_{5}}}(x),$$
for $\epsilon _{Q_{i}}\in \{1,-1\}(i=1,3,5)$.

From above, there are exactly 24 quaternary QR codes. On the other hand, for $g=2$ and Theorem 2.2, the number of all quaternary QR codes is equal to $\prod_{r=0}^{g-1}(2\times(2^{g-r}-1) )^{\dbinom{g}{r}}=24$.

Next, we consider the number of self-orthogonal or LCD quaternary QR codes. In fact, for case 1, fixed $\epsilon_{Q_{i_{r}}}\in \{1,-1\}(Q_{1}=15, Q_{3}=5, Q_{5}=3),$ we have
$$
\mathcal{C}=<g(x)>=<F_{1,15}^{\epsilon _{15}}(x)F_{3,5}^{\epsilon _{5}}(x)F_{5,3}^{\epsilon _{3}}(x)>=\prod _{j\in M^{\epsilon_{15}}_{1,15}}(x-\beta ^{j})\prod _{j\in M^{\epsilon_{5}}_{3,5}}(x-\beta ^{j})\prod _{j\in M^{\epsilon_{3}}_{5,3}}(x-\beta ^{j}).
$$
Then for each $j\in M^{\epsilon_{15}} _{1,15}\bigcup M^{\epsilon_{5}}_{3,5}\bigcup M^{\epsilon_{3}}_{5,3}, g(\beta ^{j})=0.$ By Proposition 1.1, $\mathcal{C}$ is LCD if $g(\beta ^{-j})=0,$  which means that
$$
\left(\frac{-1}{15}\right)=\left(\frac{-1}{5}\right)=\left(\frac{-1}{3}\right)=1.
$$
This is impossible for $3 \equiv -1 (\mbox{mod}\ 4).$  Similarly, from Proposition 1.2, $\mathcal{C}$ is self-orthogonal if $g(\beta ^{-j})\neq 0,$, which means that
$$
\left(\frac{-1}{15}\right)=\left(\frac{-1}{5}\right)=\left(\frac{-1}{3}\right)=-1.
$$
This is impossible for $5 \equiv 1 (\mbox{mod}\ 4).$ Hence there exist no self-orthogonal or LCD  quaternary QR codes for case 1. For cases 2-3, in the same proofs as that for case 1,
there exist no self-orthogonal or LCD quaternary QR codes $[15,8]$, which are just in accordance with the corresponding results of Theorems 2.2-2.3.

Finally, we consider the minimum distance for these 24 codes. For case 1, consider the quaternary QR code $\mathcal{C}=<g(x)>$ with
 $$
 g(x)=F_{1,15}^{1}(x)F_{3,5}^{1}(x)F_{5,3} ^{1}(x)=x^{7}+x^{6}+\alpha x^{4}+x^{2}+(\alpha +1)x+\alpha+1,
 $$
 and then the corresponding parity-check polynomial
 $$
 h(x)=(x^{n}-1)/g(x)=x^{8}+x^{7}+x^{6}+(\alpha +1)x^{5} +x^{4}+\alpha x^{3}+x^{2}+\alpha x+\alpha.
 $$
Hence we can get a parity-check matrix of $\mathcal{C}$ as follows:
\begin{equation}
H=\left(\begin{array}{ccccccccccccccc}
  1 & 1 & 1 & \alpha+1 & 1 & \alpha & 1 & \alpha & \alpha & 0 & 0 & 0 & 0 & 0 & 0  \\
  0 & 1 & 1 & 1 & \alpha+1 & 1 & \alpha & 1 & \alpha & \alpha & 0 & 0 & 0 & 0 & 0  \\
  0 & 0 & 1 & 1 & 1 & \alpha+1 & 1 & \alpha & 1 & \alpha & \alpha & 0 & 0 & 0 & 0  \\
  0 & 0 & 0 & 1 & 1 & 1 & \alpha+1 & 1 & \alpha & 1 & \alpha & \alpha & 0 & 0 & 0  \\
  0 & 0 & 0 & 0 & 1 & 1 & 1 & \alpha+1 & 1 & \alpha & 1 & \alpha & \alpha & 0 & 0  \\
  0 & 0 & 0 & 0 & 0 & 1 & 1 & 1 & \alpha+1 & 1 & \alpha & 1 & \alpha & \alpha & 0  \\
  0 & 0 & 0 & 0 & 0 & 0 & 1 & 1 & 1 & \alpha+1 & 1 & \alpha & 1 & \alpha & \alpha  \\
  \end{array}
 \right)_{7\times 15}
\end{equation}
Suppose that $b_{1},b_{2},\ldots,b_{15}$ are the columns of $H$, and then by computing, $b_{1},b_{2},b_{3},b_{4},b_{6},b_{7}$ are linearly $\mathbb{F}_4$-dependent and for any 5 columns of $H$,  they are linearly $\mathbb{F}_4$-independent, which means that the minimum distance $d=6$ of $\mathcal{C}$ by Lemma 1.1.

In the same way, one can get the minimum distance for the other 7 quaternary QR codes in case 1, which are listed in Table 2. And the minimum distance for all quaternary QR codes in cases 2-3 are listed in Tables 3-4, respectively.
\begin{table}[!h]
\centering
\caption{}
\begin{tabular}{|c|c|}
\hline
 Generator Polynomial &  Minimum Distance  \\
  $F_{1,15}^{1}(x)F_{3,5}^{1}(x)F_{5,3}^{-1}(x)=a+a x^2+x^3+x^4+a x^4+a x^5+x^7 $& 6 \\
  $F_{1,15}^{1}(x)F_{3,5}^{-1}(x)F_{5,3}^{1}(x)=1+a+x^2+a x^2+x^3+a x^4+x^5+a x^5+x^7$&  6\\
  $F_{1,15}^{1}(x)F_{3,5}^{-1}(x)F_{5,3}^{-1}(x)=a+a x+x^2+x^4+a x^4+x^6+x^7  $&6  \\
  $F_{1,15}^{-1}(x)F_{3,5}^{1}(x)F_{5,3}^{1}(x)=1+a+x^2+a x^3+x^4+a x^4+x^5+x^7     $&  6\\
  $F_{1,15}^{-1}(x)F_{3,5}^{1}(x)F_{5,3}^{-1}(x)=a+a x+x^3+a x^3+a x^5+x^6+x^7       $& 6 \\
  $F_{1,15}^{-1}(x)F_{3,5}^{-1}(x)F_{5,3}^{1}(x)=1+a+x+a x+a x^3+x^5+a x^5+x^6+x^7         $&  6 \\
  $F_{1,15}^{-1}(x)F_{3,5}^{-1}(x)F_{5,3}^{-1}(x)=a+x^2+x^3+a x^3+a x^4+x^5+x^7      $&6 \\
  $F_{1,15}^{1}(x)F_{3,5}^{1}(x)F_{5,3}^{1}(x)=1+a+x+a x+x^2+a x^4+x^6+x^7         $& 6 \\
  \hline
\end{tabular}
\end{table}
\begin{table}[!h]
\centering
\caption{}
\begin{tabular}{|c|c|}
\hline
   Generator Polynomial  &    Minimum Distance\\
  $F_{1,3}^{1}(x)F_{3,5}^{1}(x)F_{5,3} ^{-1}(x) = 1+x+a x+x^2+x^3+a x^4+x^6+a x^6+x^7 $& 6 \\
  $F_{1,3}^{1}(x)F_{3,5}^{-1}(x)F_{5,3} ^{1} (x) = a+x+a x^2+x^5+x^6+a x^6+x^7$&  3\\
  $F_{1,3}^{1}(x)F_{3,5}^{-1}(x)F_{5,3} ^{-1}(x)= 1+a x+x^3+a x^3+x^4+x^5+a x^6+x^7$ &6  \\
  $F_{1,3}^{-1}(x)F_{3,5}^{1}(x)F_{5,3} ^{1}(x)= 1+x+a x+a x^3+x^4+x^5+x^6+a x^6+x^7$ &  6\\
  $F_{1,3}^{-1}(x)F_{3,5}^{1}(x)F_{5,3} ^{-1}(x) =  1+a+x+x^2+a x^2+x^5+a x^6+x^7$ & 3 \\
  $F_{1,3}^{-1}(x)F_{3,5}^{-1}(x)F_{5,3} ^{1}(x)= 1+a x+x^2+x^3+x^4+a x^4+a x^6+x^7$ &  6 \\
  $F_{1,3}^{-1}(x)F_{3,5}^{-1}(x)F_{5,3} ^{-1}(x)= 1+a+a x+x^2+a x^2+x^5+x^6+a x^6+x^7$ &3 \\
  $F_{1,3}^{1}(x)F_{3,5}^{1}(x)F_{5,3} ^{1}(x)= a+x+a x+a x^2+x^5+a x^6+x^7$ & 3 \\
  \hline
\end{tabular}
\end{table}
\begin{table}[!h]
\centering
\caption{}
\begin{tabular}{|c|c|}
 \hline
   Generator Polynomial &  Minimum Distance \\
  $F_{1,5}^{1}(x)F_{3,5}^{1}(x)F_{5,3} ^{-1}(x)= a+x+x^3+x^4+a x^4+a x^6+x^7 $& 4 \\
  $F_{1,5}^{1}(x)F_{3,5}^{-1}(x)F_{5,3} ^{1} (x)= 1+a+a x+x^3+x^5+a x^6+x^7$&  6\\
  $F_{1,5}^{1}(x)F_{3,5}^{-1}(x)F_{5,3} ^{-1}(x) = a+x+a x+a x^2+a x^4+x^6+a x^6+x^7$ &6  \\
  $F_{1,5}^{-1}(x)F_{3,5}^{1}(x)F_{5,3} ^{1}(x) = 1+a+a x+x^2+a x^2+x^4+a x^4+a x^6+x^7$ &  6\\
  $F_{1,5}^{-1}(x)F_{3,5}^{1}(x)F_{5,3} ^{-1}(x)= a+x+a x+x^3+x^5+x^6+a x^6+x^7$ & 6 \\
  $F_{1,5}^{-1}(x)F_{3,5}^{-1}(x)F_{5,3} ^{1}(x)= 1+a+x+x^3+a x^4+x^6+a x^6+x^7$ &  4 \\
  $F_{1,5}^{-1}(x)F_{3,5}^{-1}(x)F_{5,3} ^{-1}(x)= a+x+x^3+a x^3+a x^4+a x^6+x^7$ &4 \\
  $F_{1,5}^{1}(x)F_{3,5}^{1}(x)F_{5,3} ^{1}(x)= 1+a+x+a x^3+x^4+a x^4+x^6+a x^6+x^7$ & 4 \\
  \hline
\end{tabular}
\end{table}

\end{exmp}

\begin{exmp} Consider binary QR codes with length $n=161=7\times 23$ and dimension $81$.

By Definition 2.1, the construction is described below.

First, we have $A_0=\{1\}, A_1=\{7,23\}$
$$
M_{1}=\{ i\,|\,\gcd(i,161)=1,1\leq i\leq 160\},M_{7}=\{ i\,| \,\gcd(i,n)=7, 1\leq i\leq 160\},
$$
$$
M_{23}=\{ i\,| \,\gcd(i,161)=23, 1\leq i\leq 160\},
$$
and
$$Q_{1}\in \{ 7,23,161\},Q_{7}=23 ,Q_{23}=7.$$

Case 1. By taking $Q_{1}= 161,$ we know that
$$M_{1,161}^{1}=\left\{ i\,| \,\gcd(i,161)=1, 1\leq i\leq 160, \left(\frac{i}{161}\right)=1\right\},
$$
$$
M_{1,161}^{-1}=\left\{ i\,|\gcd(i,161)=1, 1\leq i\leq 160, \left(\frac{i}{161}\right)=-1\right\},
$$
$$
M_{7,23}^{1}=\left\{ i\,| \,\gcd(i,161)=7, 1\leq i\leq 160, \left(\frac{i}{23}\right)=1\right\},
$$
$$
M_{7,23}^{-1}=\left\{ i\,| \,\gcd(i,161)=7, 1\leq i\leq 160, \left(\frac{i}{23}\right)=-1\right\},
$$
$$
M_{23,7}^{1}=\left\{ i\,| \,\gcd(i,161)=23, 1\leq i\leq 160, \left(\frac{i}{7}\right)=1\right\},
$$
$$
A_{23,7}^{-1}=\left\{ i\,| \,\gcd(i,161)=23, 1\leq i\leq 160 \left(\frac{i}{7}\right)=-1\right\},
$$
and then
$$x^{161}-1=(x-1)F_{1,161}^{1}(x)F_{1,161}^{-1}(x)F_{7,23}^{1}(x)F_{7,23}^{-1}(x)F_{23,7}^{1}(x)F_{23,7}^{-1}(x),$$
where
\begin{eqnarray*}
F_{1,161}^{1}(x)&=&(1+x+x^3+x^7+x^9+x^{11}+x^{12}+x^{13}+x^{14}+x^{17}+x^{18}+x^{19}+x^{21}+x^{23}\\
      & &+x^{24}+x^{25}+x^{26}+x^{27}+x^{28}+x^{29}+x^{33}) (1+x^2+x^3+x^5+x^8+x^9+x^{10}\\
     & &+x^{12}+x^{14}+x^{16}+x^{19}+x^{20}+x^{24}+x^{26}+x^{29}+x^{31}+x^{33}),
\end{eqnarray*}
\begin{eqnarray*}
F_{1,161}^{-1}(x)&=&(1+x^2+x^4+x^7+x^9+x^{13}+x^{14}+x^{17}+x^{19}+x^{21}+x^{23}+x^{24}+x^{25}+x^{28}+x^{30}\\
            & &+x^{31}+x^{33})(1+x^4+x^5+x^6+x^7+x^8+x^9+x^{10}+x^{12}+x^{14}+x^{15}+x^{16}+x^{19}\\
            & &+x^{20}+x^{21}+x^{22}+x^{24}+x^{26}+x^{30}+x^{32}+x^{33}),
\end{eqnarray*}
$$
F_{7,23}^{1}(x)=1+x+x^5+x^6+x^7+x^9+x^{11},\ F_{7,23}^{-1}(x)=1+x^2+x^4+x^5+x^6+x^{10}+x^{11},
$$
and
$$
F_{23,7}^{1}(x)=1+x^2+x^3,\quad F_{23,7}^{-1}(x)=1+x+x^3.
$$

Hence by Definition 2.1 we can get $8$ binary QR codes with length $161$ and dimension $81$.

Case 2. By taking $Q_{1}=7$, we know that
$$
M_{1,7}^{1}=\left\{ i\,| \,\gcd(i,161)=1, 1\leq i\leq 160,\left(\frac{i}{7}\right)=1\right\},
$$
$$
M_{1,7}^{-1}=\left\{ i\,| \,\gcd(i,161)=1, 1\leq i\leq 160 , \left(\frac{i}{7}\right)=-1\right\},
$$
and then
\begin{eqnarray*}
F_{1,7}^{1}(x)&=& (1+x+x^3+x^7+x^9+x^{11}+x^{12}+x^{13}+x^{14}+x^{17}+x^{18}+x^{19}+x^{21}+x^{23}+x^{24}\\
            & & +x^{25}+x^{26}+x^{27}+x^{28}+x^{29}+x^{33})(1+x^2+x^4+x^7+x^9+x^{13}+x^{14}+x^{17}+x^{19}\\
            & & +x^{21}+x^{23}+x^{24}+x^{25}+x^{28}+x^{30}+x^{31}+x^{33}),
\end{eqnarray*}
\begin{eqnarray*}
F_{1,7}^{-1}(x)&=&(1+x^2+x^3+x^5+x^8+x^9+x^{10}+x^{12}+x^{14}+x^{16}+x^{19}+x^{20}+x^{24}+x^{26}  \\
            & & +x^{29}+x^{31}+x^{33})(1+x^4+x^5+x^6+x^7+x^8+x^9+x^{10}+x^{12}+x^{14}+x^{15}\\
            & & +x^{16}+x^{19}+x^{20}+x^{21}+x^{22}+x^{24}+x^{26}+x^{30}+x^{32}+x^{33}).
\end{eqnarray*}
Similarly, there are exactly $8$ binary QR codes with length $161$ and dimension $81$.

Case 3. By taking  $Q_{1}=23$, we have
$$
M_{1,23}^{1}=\{ i\,| \,\gcd(i,161)=1, 1\leq i\leq 160, \left(\frac{i}{23}\right)=1\},
$$
$$
M_{1,23}^{-1}=\{ i\,| \,\gcd(i,161)=1, 1\leq i\leq 160, \left(\frac{i}{23}\right)=-1\},
$$
and then
\begin{eqnarray*}
F_{1,23}^{1}(x)&=& (1+x+x^3+x^7+x^9+x^{11}+x^{12}+x^{13}+x^{14}+x^{17}+x^{18}+x^{19}+x^{21}+x^{23}+x^{24} \\
             & & +x^{25}+x^{26}+x^{27}+x^{28}+x^{29}+x^{33})(1+x^4+x^5+x^6+x^7+x^8+x^9+x^{10}+x^{12}   \\
            & & +x^{14}+x^{15}+x^{16}+x^{19}+x^{20}+x^{21}+x^{22}+x^{24}+x^{26}+x^{30}+x^{32}+x^{33}),
\end{eqnarray*}
\begin{eqnarray*}
F_{1,23}^{-1}(x)&=& (1+x^2+x^3+x^5+x^8+x^9+x^{10}+x^{12}+x^{14}+x^{16}+x^{19}+x^{20}+x^{24} \\
             & & +x^{26}+^{29}+x^{31}+x^{33})(1+x^2+x^4+x^7+x^9+x^{13}+x^{14}+x^{17}+x^{19}\\
            & & +x^{21}+x^{23}+x^{24}+x^{25}+x^{28}+x^{30}+x^{31}+x^{33}).
\end{eqnarray*}
In the same way, we can get $8$ binary QR codes with length $161$ and dimension $81.$

From the above, there are exactly 24 binary QR codes with length $161$ and dimension $81$, which is in accordance with the corresponding result in Theorem 2.2 since
$\prod_{r=0}^{g-1}(2\times(2^{g-r}-1))^{\dbinom{g}{r}}=6\times 4=24$ when $g=2$.

Now from the proof for Proposition 1.2, we know that $\mathcal{C}$ is self-orthogonal if and only if
$$
\left(\frac{-1}{Q_{1}}\right)=\left(\frac{-1}{Q_7}\right)=\left(\frac{-1}{Q_{23}}\right)=-1,
$$
where $Q_{1}\in \{7,23,161\},Q_{7}=23$ and $Q_{23}=7.$ Note that $7\equiv 23\equiv 3(\mbox{mod}\ 4)$, so $\mathcal{C}$ is self-orthogonal if and only if $Q_{1}\in \{ 7,23\},Q_{7}=23$ and $Q_{23}=7,$ which is true only for cases 2-3. Therefore there exist 16 self-orthogonal binary QR codes, which is in accordance with the corresponding results of Theorems 2.2-2.3 since $g=2, 2^{g2^{g-1}}=16$.

Similarly, from the proof of Proposition 1.1, $\mathcal{C}$ is LCD if and only if
$$\left(\frac{-1}{Q_{1}}\right)=\left(\frac{-1}{Q_7}\right)=\left(\frac{-1}{Q_{23}}\right)=1,$$
where $Q_{1}\in \{ 7,23,161\},Q_{7}=23$ and $Q_{23}=7.$ This is impossible since $23\equiv 7\equiv 3 (\mbox{mod}\ 4).$ Hence there exists no LCD binary QR codes, which is also in accordance with the corresponding results in Theorem 2.2-2.3.
\end{exmp}

\begin{exmp}  Consider ternary QR codes with length $n=253=11\times 23$ and dimension $127$.

By Definition 2.1, the construction is described below.

First, we have $A_0=\{1\}, A_1=\{11,23\},$
$$
M_{1}=\{ i\,|\,\gcd(i,253)=1,1\leq i\leq 252\},\quad M_{11}=\{ i\,| \,\gcd(i,253)=11, 1\leq i\leq 252\},
$$
$$
M_{23}=\{ i\,| \,\gcd(i,253)=23, 1\leq i\leq 252\},
$$
and
$$
Q_{1}=\in \{ 11,23,253\}, Q_{11}=23,  Q_{23}=11.
$$

Case 1. By taking $Q_{1}=253,$
we know that
$$
M_{1,253}^{1}=\{ i\,|\gcd(i,253)=1, 1\leq i\leq 252, \left(\frac{i}{253}\right)=1\},
$$
$$
M_{1,253}^{-1}=\{ i\,|\gcd(i,253)=1, 1\leq i\leq 252, \left(\frac{i}{253}\right)=-1\},
$$
$$
M_{11,23}^{1}=\{ i|\gcd(i,253)=11, 1\leq i\leq 252, \left(\frac{i}{23}\right)=1\},
$$
$$
M_{11,23}^{-1}=\{ i| \gcd(i,253)=11, 1\leq i\leq 252, \left(\frac{i}{23}\right)=-1\},
$$
$$
M_{23,11}^{1}=\{ i\,| \,\gcd(i,253)=23, 1\leq i\leq 252, \left(\frac{i}{11}\right)=1\},
$$
$$
M_{23,11}^{-1}=\{ i\,| \,\gcd(i,253)=23, 1\leq i\leq 252 \left(\frac{i}{11}\right)=-1\},
$$
and then
$$
x^{253}-1=(x-1)F_{1,253}^{1}(x)F_{1,253}^{-1}(x)F_{11,23}^{1}(x)F_{11,23}^{-1}(x)F_{23,11}^{1}(x)F_{23,11}^{-1}(x),
$$
where
\begin{eqnarray*}
F_{1,253}^{1}(x)&=& (2+x+x^2+2x^3+x^4+x^6+2x^7+x^9+x^{10}+x^{11}+2x^{12}+2x^{16}+2x^{18}+x^{20}+x^{21} \\
             & & +x^{22}+x^{23}+x^{24}+2x^{26}+2x^{29}+x^{30}+2x^{32}+2x^{33}+x^{35}+2x^{36}+x^{37}+x^{38}+2x^{41}\\
             & & +x^{42}+2x^{43}+x^{44}+2x^{45}+x^{46}+2x^{48}+2x^{49}+x^{55}) (2+2x^2+x^3+x^4+x^5+x^6+2x^7\\
             & & +2x^8+x^9+x^{10}+2x^{12}+x^{13}+2x^{14}+2x^{15}+x^{17}+2x^{20}+2x^{21}+x^{23}+x^{24}+2x^{26}\\
             & & +2x^{27}+x^{28}+2x^{29}+x^{30}+x^{31}+2x^{32}+x^{34}+2x^{37}+x^{38}+2x^{39}+x^{41}+2x^{42}\\
             & & +x^{43}+2x^{45}+2x^{46}+x^{47}+2x^{49}+x^{50}+x^{52}+x^{55}),
\end{eqnarray*}
\begin{eqnarray*}
F_{1,253}^{-1}(x)&=&(2+2x^3+2x^5+x^6+2x^8+x^9+x^{10}+2x^{12}+x^{13}+2x^{14}+x^{16}+2x^{17}+x^{18}+2x^{21}\\
              & &+x^{23}+2x^{24}+2x^{25}+x^{26}+2x^{27}+x^{28}+x^{29}+2x^{31}+2x^{32}+x^{34}+x^{35}+2x^{38}\\
              & &+x^{40}+x^{41}+2x^{42}+x^{43}+2x^{45}+2x^{46}+x^{47}+x^{48}+2x^{49}+2x^{50}+2x^{51}+2x^{52}\\
              & &+x^{53} + x^{55}) (2+x^6+x^7+2x^9+x^{10}+2x^{11}+x^{12}+2x^{13}+x^{14}+2x^{17}+2x^{18}+x^{19}\\
              & &+2x^{20}+x^{22}+x^{23}+2x^{25}+x^{26}+x^{29}+2x^{31}+2x^{32}+2x^{33}+2x^{34}+2x^{35}+x^{37}\\
              & &+x^{39}+x^{43}+2x^{44}+2x^{45}+2x^{46}+x^{48}+2x^{49}+2x^{51}+x^{52}+2x^{53}+2x^{54}+x^{55}),
\end{eqnarray*}
$$
F_{11,23}^{1}(x)=2+x^3+x^5+2x^7+2x^8+x^9+x^{10}+x^{11},
$$
$$
 F_{11,23}^{-1}(x)=2+2x+2x^2+x^3+x^4+2x^6+2x^8+x^{11},
$$
and
$$
F_{23,11}^{1}(x)=2+x^2+2x^3+x^4+x^5,\quad F_{23,11}^{-1}(x)=2+2x+x^2+2x^3+x^5.
$$
Hence by Definition 2.1 we can get $8$ ternary QR codes with length $253$ and dimension $127$.

Case 2. By taking $Q_{1}=11$, we know that
$$
M_{1,11}^{1}=\left\{ i\,|\gcd(i,253)=1, 1\leq i\leq 252, \left(\frac{i}{11}\right)=1\right\},
$$
$$
M_{1,11}^{-1}=\left\{ i\,|\gcd(i,253)=1, 1\leq i\leq 252, \left(\frac{i}{11}\right)=-1\right\},
$$
and then
\begin{eqnarray*}
F_{1,11}^{1}(x)&=& (2+x+x^2+2x^3+x^4+x^6+2x^7+x^9+x^{10}+x^{11}+2x^{12}+2x^{16}+2x^{18}+x^{20}+x^{21}  \\
             & &  +x^{22}+x^{23} +x^{24}+2x^{26}+2x^{29}+x^{30}+2x^{32}+2x^{33}+x^{35}+2x^{36}+x^{37}\\
             & &  +x^{38}+2x^{41}+x^{42}+2x^{43} +x^{44}+2x^{45}+x^{46}+2x^{48}+2x^{49}+x^{55})  \\
             & &  (2+2x^3+2x^5+x^6+2x^8+x^9+x^{10}+2x^{12}+x^{13}+2x^{14}+x^{16}+2x^{17}+x^{18}+2x^{21} \\
             & &  +x^{23}+2x^{24}+2x^{25}+x^{26}+2x^{27}+x^{28}+x^{29}+2x^{31} +2x^{32}+x^{34}+x^{35}+2x^{38}+x^{40}\\
             & &  +x^{41}+2x^{42}+x^{43}+2x^{45}+2x^{46}+x^47+x^{48}+2x^{49}+2x^{50}+2x^{51}+2x^{52}+x^{53}+x^{55}),
\end{eqnarray*}
\begin{eqnarray*}
F_{1,11}^{-1}(x)&=& (2+2x^2+x^3+x^4+x^5+x^6+2x^7+2x^8+x^9+x^{10}+2x^{12}+x^{13}+2x^{14}+2x^{15}+x^{17} \\
              & & +2x^{20}+2x^{21}+x^{23}+x^{24}+2x^{26}+2x^{27}+x^{28}+2x^{29}+x^{30}+x^{31}+2x^{32}+x^{34}+2x^{37}\\
              & & +x^{38}+2x^{39}+x^{41}+2x^{42}+x^{43}+2x^{45}+2x^{46}+x^{47}+2x^{49}+x^{50}+x^{52}+x^{55})\\
             & &  (2+x^6+x^7+2x^9+x^{10}+2x^{11}+x^{12}+2x^{13}+x^{14}+2x^{17} +2x^{18}+x^{19}+2x^{20}+x^{22}\\
              & &  +x^{23}+2x^{25}+x^{26}+x^{29}+2x^{31}+2x^{32}+2x^{33}+2x^{34}+2x^{35}+x^{37}+x^{39}+x^{43}\\
             & &  +2x^{44}+2x^{45}+2x^{46}+x^{48}+2x^{49}+2x^{51}+x^{52}+2x^{53}+2x^{54}+ x^{55}),
\end{eqnarray*}
Similarly, there are exactly $8$ ternary QR codes with length $253$ and dimension $127$.

Case 3. By taking $Q_{1}=23$, we have
$$
M_{1,23}^{1}=\left\{ i\,|\gcd(i,253)=1, 1\leq i\leq 252, \left(\frac{i}{23}\right)=1\right\},
$$
$$
M_{1,23}^{-1}=\left\{ i\,|\gcd(i,253)=1, 1\leq i\leq 252, \left(\frac{i}{23}\right)=-1\right\}$$
and then
\begin{eqnarray*}
F_{1,23}^{1}(x)&=&  2+x+x^2+2x^3+x^4+x^6+2x^7+x^9+x^{10}+x^{11}+2x^{12}+2x^{16}+2x^{18}+x^{20}+x^{21} \\
              & & +x^{22}+x^{23}+x^{24}+2x^{26}+2x^{29}+x^{30}+2x^{32}+2x^{33}+x^{35}+2x^{36}+x^{37}+x^{38} \\
              & & +2x^{41}+x^{42}+2x^{43}+x^{44}+2x^{45}+x^{46}+2x^{48}+2x^{49}+x^{55}) (2+x^6+x^7+2x^9+x^{10}  \\
             & &  +2x^{11}+x^{12}+2x^{13}+x^{14}+2x^{17}+2x^{18}+x^{19}+2x^{20}+x^{22}+x^{23}+2x^{25}+x^{26}  \\
              & & +x^{29}+2x^{31}+2x^{32}+2x^{33}+2x^{34}+2x^{35}+x^{37}+x^{39}+x^{43}+2x^{44}+2x^{45}+2x^{46}\\
             & &  +x^{48}+2x^{49}+2x^{51}+x^{52}+2x^{53}+2x^{54}+x^{55}),
\end{eqnarray*}
\begin{eqnarray*}
F_{1,23}^{-1}(x)&=& (2+2x^3+2x^5+x^6+2x^8+x^9+x^{10}+2x^{12}+x^{13}+2x^{14}+x^{16}+2x^{17}+x^{18}+2x^{21}   \\
              & & +x^{23}+2x^{24}+2x^{25}+x^{26}+2x^{27}+x^{28}+x^{29}+2x^{31}+2x^{32}+x^{34}+x^{35}\\
              & & +2x^{38}+x^{40}+x^{41}+2x^{42}+x^{43}+2x^{45}+2x^{46}+x^{47}+x^{48}+2x^{49}+2x^{50}\\
              & & +2x^{51}+2x^{52}+x^{53}+x^{55})(2+2x^2+x^3+x^4+x^5+x^6+2x^7+2x^8+x^9+x^{10}+2x^{12}\\
             & &  +x^{13}+2x^{14}+2x^{15}+x^{17}+2x^{20}+2x^{21}+x^{23}+x^{24}+2x^{26}+2x^{27}+x^{28}  \\
              & & +2x^{29}+x^{30}+x^{31}+2x^{32}+x^{34}+2x^{37}+x^{38}+2x^{39}+x^{41}+2x^{42}+x^{43}\\
             & &  +2x^{45}+2x^{46}+x^{47}+2x^{49}+x^{50}+x^{52}+x^{55})  ,
\end{eqnarray*}

In the same way, we can get $8$ ternary QR codes with length of $253$ and dimension $127$.

From the above, there are exactly 24 ternary QR codes with length $253$ and dimension $127,$
which is in accordance with the corresponding result in Theorem 2.2 since $\prod_{r=0}^{g-1}(2\times(2^{g-r}-1) )^{\dbinom{g}{r}}=6\times 4=24$ when $g=2$.

Now from the proof for Proposition 1.2, we know that $\mathcal{C}$ is self-orthogonal if and only if
$$\left(\frac{-1}{Q_{1}}\right)=\left(\frac{-1}{Q_{11}}\right)=\left(\frac{-1}{Q_{23}}\right)=-1,$$
where $Q_{1}\in \{ 11,23,253\},Q_{11}=23$ and $Q_{23}=11.$ Note that $11\equiv 23\equiv 3(\mod 4),$ so $\mathcal{C}$ is self-orthogonal if and only if  $Q_{1}\in \{ 11,23\},Q_{11}=23$ and $Q_{23}=11,$ which is true only for cases 2-3. Therefore there exist 16 self-orthogonal ternary QR codes, which is in accordance with the corresponding results of Theorem 2.2-2.3 since $g=2$, $2^{g2^{g-1}}=16$.

Similarly, from the proof of  Proposition 1.1, $\mathcal{C}$ is LCD if and only if
$$\left(\frac{-1}{Q_{1}}\right)=\left(\frac{-1}{Q_{11}}\right)=\left(\frac{-1}{Q_{23}}\right)=1,$$
where $Q_{1}\in \{ 11,23,253\},Q_{11}=23$ and $Q_{23}=11.$ This is impossible since $23\equiv 11\equiv 3 (\mod 4).$ And so there exists no LCD ternary QR codes, which also is in accordance with the corresponding results in Theorem 2.2-2.3.
\end{exmp}

\begin{exmp}
\quad By taking $p_1=5,p_2=17,p_3=29,n=p_1p_2p_3=2465$, and $\theta$ a 2465-th primitive root of the unity in $\mathbb{F}_{2^{56}}$. Then by Theorem 2.1, there exists a $q$-ary QR code if and only if $\left(\frac{q}{5}\right)=\left(\frac{q}{17}\right)=\left(\frac{q}{29}\right)=1$. In particular, it is obvious for $q=4$, therefore there is a quaternary QR code $[2465, 1233]$.

Similarly, there exists a LCD QR code with length 2465 and dimension $\frac{n+1}{2}=1233$ if and only if for some
$$
Q_1\in\{5,17,29,85,145,493,2465\},\ Q_{5}\in\{17,29,493\},\ Q_{17}\in\{5,29,145\},\ Q_{29}\in\{5,17,85\},
$$
we have
$$
\left(\frac{-1}{Q_1}\right)=\left(\frac{-1}{Q_{5}}\right)=\left(\frac{-1}{Q_{17}}\right)
=\left(\frac{-1}{Q_{29}}\right)=\left(\frac{-1}{5}\right)=\left(\frac{-1}{17}\right)=\left(\frac{-1}{29}\right)=1,
$$
which is always true since $5\equiv 17\equiv29\equiv 1(\mbox{mod}\ 4)$. This means that any QR code is LCD and there exists no self-orthogonal QR code with length 2465 and dimension 1233.
\end{exmp}

\begin{exmp}
\quad By taking $p_1=3,p_2=7,p_3=11,n=p_1p_2p_3=231$, and $\theta$ a 231-th primitive root of the unity in $\mathbb{F}_{2^{30}}$. Then by Theorem 2.1, there exists a $q$-ary QR code with length 231 and dimension
116 if and only if $\left(\frac{q}{3}\right)=\left(\frac{q}{7}\right)=\left(\frac{q}{11}\right)=1$. In particular, it is obvious for $q=4$, therefore there is a quaternary QR code $[231,116]$.

Similarly, there is no LCD QR code with length 231 and dimension 116 since
$$
\left(\frac{-1}{3}\right)=\left(\frac{-1}{7}\right)=\left(\frac{-1}{11}\right)=-1.
$$
And there exists a self-orthogonal QR code with length 231 and dimension
$\frac{n+1}{2}=116$ if and only if for some
$$
Q_1\in\{3,7,11,21,33,77,231\},\ Q_{3}\in\{7,11,77\},\ Q_{7}\in\{3,11,33\},\ Q_{11}\in\{3,7,21\},
$$
we have
$$
\left(\frac{-1}{Q_1}\right)=\left(\frac{-1}{Q_{3}}\right)=\left(\frac{-1}{Q_{7}}\right)
=\left(\frac{-1}{Q_{11}}\right)=\left(\frac{-1}{3}\right)=\left(\frac{-1}{7}\right)=\left(\frac{-1}{11}\right)=-1,
$$
which is true if and only if
$$
Q_1\in\{3,7,11,231\}, \ Q_{3}\in \{7,11\}, \ Q_{7}\in\{3,11\},\ Q_{11}\in\{3,7\}.
$$
 This means that the number of self-orthogonal QR codes with length 231 and dimension 116 is equal to
 $$
(2\times 4)\times (2\times 2)^{3} \times 2^{3}=2^{12}=2^{g2^{g-1}},
 $$
where the last equality is due to the assumption $g=3$.
\\
\end{exmp}

\section{$m$-th Residue Codes With an Odd Prime Length}

For a positive integer $m\geq 2$, this section generalizes the quadratic residue code with an odd prime length over the finite field $\mathbb{F}_{q}$ to be the case $m$-th residue code, and gives a criterion for LCD or self-orthogonal $m$-th residue codes with an odd prime length and the corresponding enumerators.

\subsection{Construction of $m$-th Residue Codes}
Let $n=p$ be an odd prime, $m\geq 2$ be a positive integer with $m|p-1$, and $q$ be a power of the prime with $p|q-1$. Suppose that $\theta$ is a $p$-th primitive root of  unity in $\mathbb{F}_q$ and $r$ is a positive integer with $ord_pr=p-1$ and $1\leq r\leq p-1$, where $ord_pr=\min\{x\,|\, r^x\equiv 1(\mbox{mod}\ p),x\in \mathbb{Z}^+ \}.$

First, set
$$
A_0=\{k\,|\,  k^{\frac{p-1}{m}}\equiv 1(\mbox{mod}\ p), 1\leq k\leq p-1\},
$$
and then $A_0$ is a cyclic subgroup of the multiplication group $\mathbb{F}^*_p$ with $|A_0|=\frac{p-1}{m}$. So there are exactly $m$ cosets $A_i(0\leq i\leq m-1)$ of $A_0$ in $\mathbb{F}^*_p$, namely,
$$
A_i=r^iA_0=\{x\,|\, x=r^ik(\mbox{mod}\ p),  k\in A_0,1\leq x\leq p-1\}.
$$
Thus we can get the disjoint union $\mathbb{F}^*_{p}=\cup^{m-1}_{i=0}A_i,$ and then the factorization
$$
x^p-1=(x-1)\prod^{m-1}_{j=0}\prod_{i\in A_j}(x-\theta^i).
$$

Now for any $j=0,1,\ldots,m-1$, set $f_j(x)=\prod_{i\in A_j}(x-\theta^i)$, then $\mathcal{C}=<f_j(x)>$ or $<(x-1)f_j(x)>$ is an ideal of the ring $\mathbb{F}_q/<x^p-1>$. Thus we have the following

\textbf{Definition 3.1}\quad  $\mathcal{C}=<f_j(x)>$ or $<(x-1)f_j(x)>$ is an $m$-th residue code with length $p$.

It is easy to see that for any $j=0,1,\ldots, m-1,\deg(f_j(x))=\frac{p-1}{m}$, and then any $m$-th residue code with length $p$ has the dimension
$$
p-\frac{p-1}{m}, \ \mbox{or}\quad p-1-\frac{p-1}{m},
$$
respectively. And there are exactly $m$  $m$-th residue codes with length $p$.

\textbf{Remark 3.2}\quad By taking $m=2$ in Definition 3.1, one can get 2 QR codes with length $p$ and the maximum dimension $\frac{p+1}{2}$, which is the corresponding result given by Prange in 1958.

\subsection{Main Results}

\textbf{Theorem 3.3}\quad (1)\quad There exists an LCD $m$-th residue code with length $p$ if and only if
$$
p\equiv 1\ (\mbox{mod}\ 2m).
$$

(2)\quad The number of LCD  $m$-th residue codes with length $p$ is $0$ or $m$.

(3)\quad There exists a self-orthogonal $m$-th residue code with length $p$ if and only if
$$
p\equiv m+1 \ (\mbox{mod}\ 2m).
$$

(4)\quad The number of $m$-th self-orthogonal residue codes with length $p$ is $0$ or $m$.

\textbf{Proof}\quad Note that $\theta$ is a $p$-th primitive root of unity, thus for $j=0,1,\ldots,m-1$, any root of $f_j(x)$ is in the form $\beta=\alpha^t$ for some $t=1,\ldots,p-1$.

(1)\quad Suppose that $\mathcal{C}=<f_j(x)>$ is LCD, then by Proposition 1.1 we know that for any $\beta, f_j(\beta)=0$ implies that $f_j(\beta^{-1})=0$, equivalently, $i\in A_j$ means that $-i\in A_j$. Therefore from the definition of $A_j$, there exists an LCD $m$-th residue code with length $p$ if and only if $(-1)^{\frac{p-1}{m}}\equiv 1(\mbox{mod}\ p),$
namely, $p\equiv 1(\mbox{mod}\ 2m)$. Furthermore, in this case, any $m$-th residue code is LCD.

This completes the proof of (1).

(2)\quad From the proof of (1), (2) is immediate.

(3)\quad If $\mathcal{C}=<f_j(x)>$ is self-orthogonal, then by Proposition 1.2 we know that for any $\beta, f_j(\beta)=0$ implies that $f_j(\beta^{-1})\neq 0$, equivalently, $i\in A_j$ means that $-i\not\in A_j$. Therefore from the definition of $A_j$, there exists a self-orthogonal $m$-th residue code with length $p$ if and only if
$(-1)^{\frac{p-1}{m}}\not\equiv 1(\mbox{mod}\ p),$ which means that $\frac{p-1}{m}\equiv 1(\mbox{mod}\ 2)$, namely, $p\equiv m+1(\mbox{mod}\ 2m)$. Furthermore, in this case, any $m$-th residue code is self-orthogonal.

This completes the proof of (3).

(4)\quad From the proof of (3), (4) is immediate.
\hfill$\Box$\\

\subsection{Example}

In this section, two examples are given to support Theorems 3.3. We obtain all 3-rd residue codes with length 7 and 4-th residue codes with length 17, respectively.

\textbf{Example 3.1}\quad Let $m=3$ and $n=7$, then $n\equiv 1(\mbox{mod}\ 2m)$. Suppose that $q$ is a power of the prime and $7|q-1, \theta$ is a fixed 7-th primitive root of unity in $\mathbb{F}_q$. By computing directly, we have $\mbox{ord}_73=6$.

First, it is easy to see that the set of 3rd residues modulo 7 is
$$
A_0=\{x\,|\, x^{\frac{7-1}{3}}=x^2\equiv 1(\mbox{mod}\ 7),1\leq x\leq 6\}=\{1,6\},
$$
and then there are exactly two cosets of $A_0$ in $\mathbb{F}^*_7$, namely,
$$
A_1=3A_0=\{3,4\},\quad \mbox{and}\quad A_2=3^2A_0=\{2,5\}.
$$
Thus we have the disjoint union $\{1,2,3,4,5,6\}=A_0\cup A_1\cup A_2,$ and then the factorization
$$
x^7-1=(x-1)\prod_{i\in A_0}(x-\theta^i)\prod_{i\in A_1}(x-\theta^i)\prod_{i\in A_2}(x-\theta^i).
$$
Therefore there are three 3rd residue codes with length 7 as follows,
$$
\mathcal{C}_j=<f_j(x)=\prod_{i\in A_j}(x-\theta^i)>(j=0,1,2),
$$
where
$$
f_0(x)=(x-\theta)(x-\theta^6), \ f_1(x)=(x-\theta^3)(x-\theta^4),\ f_0(x)=(x-\theta^2)(x-\theta^5).
$$

 It is easy to see that for any $j=0,1,2, f_j(\theta^t)=0$ if and only if $f_j(\theta^{7-t})=0$. Note that $\theta$ is a 7-th primitive root of unity, namely, $(\theta^t)^{-1}=\theta^{-t}=\theta^{7-t}$ for any $t\in\{0,1,\ldots,6\}$. Now by Propositions 1.1-1.2, we know that any $\mathcal{C}_j(0\leq j\leq 2)$ is LCD but not self-orthogonal.

\textbf{Example 3.2}\quad Let $m=4$ and $n=17$, then $n\equiv 1(\mbox{mod}\ 2m)$. Suppose that $q$ is a power of the prime with $n|q-1$, and $\theta$ is a 17-th primitive root of unity in $\mathbb{F}_q$. By directly computing, we have $\mbox{ord}_{17}3=16$.

First, it is easy to see that the set of 4-th residues modulo 17 is
$$
A_0=\{x\,|\, x^{\frac{17-1}{4}}=x^4\equiv 1(\mbox{mod}\ 17),1\leq x\leq 16\}=\{1,4,13,16\},
$$
and then there are exactly 3 cosets of $A_0$ modulo 17, namely,
$$
A_1=3A_0=\{3,5,12,14\},\quad A_2=3^2A_0=\{2,8,9,15\}, \quad \mbox{and}\quad A_3=3^3A_0=\{6,7,10,11\},
$$
Thus we have the disjoint union
$$
\{1,2,\ldots,16\}=A_0\cup A_1\cup A_2\cup A_3,
$$
and then the factorization
$$
x^{17}-1=(x-1)\prod_{i\in A_0}(x-\theta^i)\prod_{i\in A_1}(x-\theta^i)\prod_{i\in A_2}(x-\theta^i)\prod_{i\in A_3}(x-\theta^i).
$$
Therefore there are four 4-th residue codes with length 17 as follows,
$$
\mathcal{C}_j=<f_j(x)=\prod_{i\in A_j}(x-\theta^i)>(j=0,1,2,3),
$$
where
$$
f_0(x)=(x-\theta)(x-\theta^4)(x-\theta^{13})(x-\theta^{16}), \ f_1(x)=(x-\theta^3)(x-\theta^5)(x-\theta^{12})(x-\theta^{14}),
$$
$$
f_2(x)=(x-\theta^2)(x-\theta^8)(x-\theta^9)(x-\theta^{15}),\ f_3(x)=(x-\theta^6)(x-\theta^7)(x-\theta^{10})(x-\theta^{11}).
$$

It is easy to see that for any $j=0,1,2,3, f_{j}(\theta^t)=0$ if and only if $f_{j}(\theta^{-t})=0$, hence by Propositions 1.1-1.2, any $\mathcal{C}_j(0\leq j\leq 3)$ is LCD but not self-orthogonal.\\

\section{Conclusions and Acknowledgment}

In this paper,  for a given positive integer $m\geq 2$ and distinct odd primes $p, p_{1},p_{2},\ldots,p_{g}$, we generalize the methods in \cite{C.P,D}, and construct quadratic residue codes with length $p_{1}p_{2}\cdots p_{g}$ or $m$-th residue codes with length $p$, respectively. We also give a sufficient and necessary condition for the existence of these self-orthogonal or LCD codes and then obtain the corresponding counting formula.

The authors would sincerely like to thanks Professor Yaotsu Chang and Professor Chongdao Lee for some beneficial discussions and help compute the minimum distance in Example 2.1.


\end{document}